\documentclass[12pt]{amsart}

\usepackage{fullpage}
\usepackage{setspace}
\usepackage[leqno]{amsmath}
\usepackage{amsfonts}
\usepackage{amssymb}
\usepackage{amsthm}
\usepackage{amssymb}
\usepackage{bbm}
\usepackage{color}

\theoremstyle{plain}
\newtheorem{theorem}{Theorem}[section]

\theoremstyle{definition}

\newcommand{\ds}{\displaystyle\sum}
\newcommand{\dss}{\displaystyle}
\newcommand{\R}{\mathbb{R}}
\newcommand{\C}{\mathbb{C}}
\newcommand{\ga}{\gamma}
\newcommand{\de}{\delta}
\newcommand{\la}{\lambda}
\newcommand{\be}{\beta}

\newcommand{\si}{\sigma}
\newcommand{\ta}{\tau}
\newcommand{\ro}{\rho}
\newcommand{\al}{\alpha}
\newcommand{\Aa}{A_\alpha}
\newcommand{\Ca}{C_\alpha}
\newcommand{\Da}{D_\alpha}
\newcommand{\Ea}{E_\alpha}
\newcommand{\Fa}{F_\alpha}
\newcommand{\Ga}{G_\alpha}
\newcommand{\Ma}{M_\alpha}
\newcommand{\Na}{N_\alpha}
\newcommand{\Xa}{X_\alpha}

\newcommand{\Xb}{X_\beta}

\newcommand{\Sa}{S_\alpha}
\newcommand{\Sb}{S_\beta}
\DeclareMathOperator{\nr}{nilrad}
\DeclareMathOperator{\Ann}{Ann}

\title{Solvable Leibniz Algebras with Heisenberg Nilradical}

\author[Bosko-Dunbar]{Lindsey Bosko-Dunbar}
\address{Department of Mathematics, Spring Hill College\\
Mobile, AL 36608}
\email{lboskodunbar@shc.edu}

\author[Dunbar]{Jonathan D. Dunbar}
\address{Department of Mathematics, Spring Hill College\\
Mobile, AL 36608}
\email{jdunbar@shc.edu}

\author[Hird]{J.T. Hird}
\address{Department of Mathematics, West Virginia University, Institute of Technology\\
Montgomery, WV 25136}
\email{John.Hird@mail.wvu.edu}

\author[Stagg]{Kristen Stagg}
\address{Department of Mathematics, The University of Texas at Tyler\\
Tyler, TX 75799}
\email{kstagg@uttyler.edu}

\begin{document}

\subjclass[2010]{17D99}
\keywords{Leibniz, Heisenberg, nilradical, classification, Lie}

\doublespacing
\maketitle

\begin{abstract}
All solvable Lie algebras with Heisenberg nilradical have already been classified.  We extend this result to a classification of solvable Leibniz algebras with Heisenberg nilradical.  As an example, we show the complete classification of all real or complex Leibniz algebras whose nilradical is the 3-dimensional Heisenberg algebra.
\end{abstract}

\section{Introduction}\label{intro}


Leibniz algebras were defined by Loday in 1993 \cite{loday, loday2}.  In recent years it has been a common theme to extend various results from Lie algebras to Leibniz algebras \cite{ao, ayupov, omirov}.  Several authors have proven results on nilpotency and related concepts which can be used to help extend properties of Lie algebras to Leibniz algebras.
Specifically, variations of Engel's theorem for Leibniz algebras have been proven by different authors \cite{barnesengel, jacobsonleib} and Barnes has proven Levi's theorem for Leibniz algebras \cite{barneslevi}.  Additionally, Barnes has shown that left-multiplication by any minimal ideal of a Leibniz algebra is either zero or anticommutative \cite{barnesleib}.

In an effort to classify Lie algebras, many authors place various restrictions on the nilradical \cite{cs, nw, tw, wld}.  In \cite{rw}, Rubin and Winternitz study solvable Lie algebras with Heisenberg nilradical.  It is the goal of this paper to extend these results to the Leibniz setting.

Recent work has been done on classification of certain classes of Leibniz algebras \cite{aor, chelsie-allison, clok, clok2}.  Especially useful for this paper is the work by Ca\~nete and Khudoyberdiyev \cite{ck} which classifies all non-nilpotent 4-dimensional Leibniz algebras over $\C$.  We recover their results for the 4-dimensional solvable Leibniz algebra over $\C$ with Heisenberg nilradical, and extend this case to classify such Leibniz algebras over $\R$.  Our result is not limited to only the 4-dimensional case, and can be used to classify all Leibniz algebras whose nilradical is Heisenberg.  In the last section we classify all (indecomposable) solvable Leibniz algebras whose nilradical is the 3-dimensional Heisenberg algebra, $H(1)$.

\section{Preliminaries}

A Leibniz algebra, $L$, is a vector space over a field (which we will take to be $\C$ or $\R$) with a bilinear operation (which we will call multiplication) defined by $[x,y]$ which satisfies the Leibniz identity
\begin{equation}\label{Jacobi}
[x,[y,z]] = [[x,y],z] + [y,[x,z]]
\end{equation}
for all $x,y,z \in L$.  In other words $L_x$, left-multiplication by $x$, is a derivation.  Some authors choose to impose this property on $R_x$, right-multiplication by $x$, instead.  Such an algebra is called a ``right'' Leibniz algebra, but we will consider only ``left'' Leibniz algebras (which satisfy \eqref{Jacobi}).  $L$ is a Lie algebra if additionally $[x,y]=-[y,x]$.

The derived series of a Leibniz (Lie) algebra $L$ is defined by $L^{(1)}=[L,L]$, $L^{(n+1)}=[L^{(n)},L^{(n)}]$ for $n\ge 1$.  $L$ is called solvable if $L^{(n)}=0$ for some $n$. The lower-central series of $L$ is defined by $L^2 = [L,L]$, $L^{n+1}=[L,L^n]$ for $n>1$. $L$ is called nilpotent if $L^n=0$ for some $n$.  It should be noted that if $L$ is nilpotent, then $L$ must be solvable.

The nilradical of $L$ is defined to be the (unique) maximal nilpotent ideal of $L$, denoted by $\nr(L)$.  It is a classical result that if $L$ is solvable, then $L^2 = [L,L] \subseteq \nr(L)$.  From \cite{mubar}, we have that
\begin{equation*}\label{dimension}
\dim (\nr(L)) \geq \frac{1}{2} \dim (L).
\end{equation*}

The Heisenberg algebra $H(n)$ is the $(2n+1)$-dimensional Lie algebra with basis \newline $\{H, P_1, \ldots, P_n, B_1, \ldots, B_n\}$ defined by multiplications
\begin{equation}\label{Heis}
[P_i,B_j]=-[B_j,P_i]=\de_{ij}H
\end{equation}
with all other products equal to zero.

The left-annihilator of a Leibniz algebra $L$ is the ideal $\Ann_\ell(L) = \left\{x\in L\mid [x,y]=0\ \forall y\in L\right\}$. Note that the elements $[x,x]$ and $[x,y] + [y,x]$ are in $\Ann_\ell(L)$, for all $x,y\in L$, because of \eqref{Jacobi}.

An element $x$ in a Leibniz algebra $L$ is nilpotent if both $(L_x)^n = (R_x)^n = 0$ for some $n$.  In other words, for all $y$ in $L$
\begin{equation*}
[x,\cdots[x,[x,y]]] = 0 = [[[y,x],x]\cdots,x].
\end{equation*}

A set of matrices $\{\Xa\}$ is called linearly nilindependent if no non-zero linear combination of them is nilpotent.  In other words, if
\begin{equation*}
X = \displaystyle\sum_{\al=1}^f c_\al \Xa,
\end{equation*}
then $X^n=0$ implies that $c_\al=0$ for all $\al$.  A set of elements of a Leibniz algebra $L$ is called linearly nilindependent if no non-zero linear combination of them is a nilpotent element of $L$.

\section{Classification}

Let $H(n)$ be the $(2n+1)$-dimensional Heisenberg (Lie) algebra over the field $F$ ($\C$ or $\R$) with basis $\{H, P_1, \ldots, P_n, B_1, \ldots, B_n\}$ and products given by \eqref{Heis}.  We will extend $H(n)$ to a solvable Leibniz algebra of dimension $2n+1+f$ by appending linearly nilindependent elements $\{S_1, \ldots, S_f\}$.  In doing so, we will construct an indecomposable Leibniz algebra whose nilradical is $H(n)$.

The left and right actions of each $\Sa$ ($\al = 1, \ldots , f$) on H(n) are given by the following matrices,
\begin{eqnarray*}\label{Ls}
L_{\Sa} \left(
\begin{array}{c}
H\\
P\\
B
\end{array}
\right) &=& \left(
\begin{array}{ccc}
2a_\al	&	\si_1^\al	&	\si_2^\al	\\
\ga_1^\al	&	a_\al I_n+\Aa	&	\Ca	\\
\ga_2^\al	&	\Da	&	a_\al I_n+\Ea
\end{array}
\right) \left(
\begin{array}{c}
H\\
P\\
B
\end{array}
\right)\\
\label{Rs}
R_{\Sa} \left(
\begin{array}{c}
H\\
P\\
B
\end{array}
\right) &=& \left(
\begin{array}{ccc}
2b_\al\	&	\ta_1^\al	&	\ta_2^\al	\\
\ro_1^\al	&	b_\al I_n+\Fa	&	\Ga	\\
\ro_2^\al	&	\Ma	&	b_\al I_n+\Na 
\end{array}
\right) \left(
\begin{array}{c}
H\\
P\\
B
\end{array}
\right)
\end{eqnarray*}
where $a_\al, b_\al \in F$;\quad $\si_i^\al, \ta_i^\al \in F^{1 \times n}$ and $\ga_i^\al, \ro_i^\al \in F^{n \times 1}$ for $i=1$ or 2;\quad $I_n$ is the $n \times n$ identity matrix;\quad  $\Aa, \Ca, \Da, \Ea, \Fa, \Ga, \Ma, \Na \in F^{n \times n}$;\quad $P = (P_1, \ldots, P_n)^T$ and $B = (B_1, \ldots, B_n)^T$.  Here, $T$ represents the transpose of the matrix.  Since $[L,L] \subseteq \nr(L)$, we have that
\begin{equation}\label{SaSb}
[\Sa, \Sb] = r_{\al\be} H + \ds_{i=1}^n \mu_{\al\be}^i P_i + \ds_{i=1}^n \nu_{\al\be}^i B_i
\end{equation}
where $r_{\al\be}, \mu_{\al\be}^i, \nu_{\al\be}^i \in F$.

By performing a change of basis to $\tilde{S}_\al = \Sa + \ds_{i=1}^n {(\ga_1^\al)}_i B_i - \ds_{i=1}^n {(\ga_2^\al)}_i P_i$, we eliminate the $H$ component of $[\Sa, P_i]$ and $[\Sa, B_i]$ in the new basis, so $\ga_1^\al = \ga_2^\al = 0$ for each $\al = 1,\ldots,f$. The Leibniz identity on the following triples imposes further constraints on $\dss L_{\Sa}$ and $\dss R_{\Sa}$.
\begin{align}\label{identities}
\begin{tabular}{ccc}
Leibniz Identity	&	&	Constraint	\\
\hline
$\{\Sa, P_i, H\}$	&	$\Rightarrow$	&	$\si_2^\al = 0$	\\
$\{\Sa, B_i, H\}$	&	$\Rightarrow$	&	$\si_1^\al = 0$	\\
$\{P_i, H, \Sa\}$	&	$\Rightarrow$	&	$\ta_2^\al = 0$	\\
$\{B_i, H, \Sa\}$	&	$\Rightarrow$	&	$\ta_1^\al = 0$	\\
$\{\Sa, P_i, P_j\}$	&	$\Rightarrow$	&	$\Ca = \Ca^T$	\\
$\{\Sa, B_i, B_j\}$	&	$\Rightarrow$	&	$\Da = \Da^T$	\\
$\{\Sa, B_i, P_j\}$	&	$\Rightarrow$	&	$\Ea = -\Aa^T$ \\
$\{\Sa, \Sb, P_i\}$	&	$\Rightarrow$	&	$\nu_{\al\be}^i = 0$	\\
$\{\Sa, \Sb, B_i\}$	&	$\Rightarrow$	&	$\mu_{\al\be}^i = 0$
\end{tabular}
\end{align}
For instance, the Leibniz identity on the triple $\{\Sa, P_i, H\}$ gives the following.
\begin{eqnarray*}
[\Sa,[P_i,H]] &=& [[\Sa,P_i],H]+[P_i,[\Sa,H]]\\
0 = [\Sa,0] &=& \left[(\ga_1^\al)_i H + a_\al P_i + \ds_{j=1}^n A_{ij}^\al P_j + \ds_{j=1}^n C_{ij}^\al B_j \ , \ H\right]\\
&& \phantom{0} + \left[P_i \ , \ 2a_\al H + \ds_{j=1}^n (\si_1^\al)_j P_j + \ds_{j=1}^n (\si_2^\al)_j B_j\right]\\
0 &=& 0 + \ds_{j=1}^n (\si_2^\al)_j [P_i, B_j]\\
0 &=& (\si_2^\al)_i H\\
0 &=& (\si_2^\al)_i
\end{eqnarray*}
Since this holds for all $i = 1, \ldots, n$, we have that $\si_2^\al = 0$.

As a result of \eqref{identities}, we observe that \eqref{SaSb} simplifies to $[\Sa,\Sb] = r_{\al\be}H$, and the matrix $\Xa = \left(\begin{array}{cc} \Aa&\Ca\\ \Da&-\Aa^T\end{array}\right)$ is an element of the symplectic Lie algebra $sp(2n,F)$.  This implies that $\Xa K + K \Xa^T = 0$, where $K = \left(\begin{array}{cc} 0 & I_n\\ -I_n & 0 \end{array}\right)$. 

We could study the components of $R_{\Sa}$, right-multiplication by $\Sa$, in the same way, but it is more convenient for our purposes to instead consider the left-annihilator $\Ann_\ell(L)$.  Since $[x,x]$,$[x,y]+[y,x]\in \Ann_\ell(L)$, for all $x,y\in L$, we get the following relations between the components of $R_{\Sa}$ and the components of $L_{\Sa}$.
\begin{align*}
\begin{tabular}{ccc}
Identity	&	&	Constraints	\\
\hline
$[[\Sa,P_i]+[P_i,\Sa]\, ,P_j] = 0$	&	$\Rightarrow$	&	$\Ga = -\Ca$\\
$[[\Sa,B_i]+[B_i,\Sa]\, ,B_j] = 0$	&	$\Rightarrow$	&	$\Ma = -\Da$\\
$\left.\begin{tabular}{l}
$[[\Sa,P_i]+[P_i,\Sa]\, , B_j] = 0$\\
$[[\Sa,B_i]+[B_i,\Sa]\, , P_j] = 0$
\end{tabular}\right\}$				&	$\Rightarrow$	&	$a_\al = -b_\al\,$, $\Fa = -\Aa$
\end{tabular}
\end{align*}
Now, defining $\ro^\al = \left(\ro_1^\al\,, \ro_2^\al\right)^T$ and using the aforementioned $\Xa$, our matrices $\dss L_{\Sa}$ and $\dss R_{\Sa}$ take the following forms.
\begin{eqnarray*}
L_{\Sa} &=& \left(\begin{array}{cc} \phantom{-}2a_\al&0\\0&\phantom{-}a_\al I_{2n}+\Xa \end{array}\right)	\\
R_{\Sa} &=& \left(\begin{array}{cc} -2a_\al&0\\ \ro^\al&-a_\al I_{2n}-\Xa \end{array}\right)
\end{eqnarray*}

Since $[\Sa,\Sb]$ is in the span of $H$ and $H$ is central, the Leibniz identity on $\{\Sa, \Sb, Z\}$ (where $Z$ is an arbitrary element of H(n)) gives that $L_{\Sa} L_{\Sb} = L_{\Sb} L_{\Sa}$.  This implies that $\Xa \Xb = \Xb \Xa$, so the $\Xa$'s commute.  Similarly, the Leibniz identity on $\{\Sa, Z, \Sb\}$ (where $Z$ is an arbitrary element of H(n)) gives that $L_{\Sa} R_{\Sb} = R_{\Sb} L_{\Sa}$.  By looking at the matrix representation of this identity (with $\al = \be$), we obtain the following.  For simplicity, we suppress the $\al$ notation and consider only one fixed $\Sa$, $S$.
\begin{eqnarray*}
\left(\begin{array}{cc} 2a & 0 \\ 0 & a I + X \end{array}\right)
\left(\begin{array}{cc} -2a & 0 \\ \ro & -a I - X \end{array}\right) &=&
\left(\begin{array}{cc} -2a & 0 \\ \ro & -a I - X \end{array}\right)
\left(\begin{array}{cc} 2a & 0 \\ 0 & a I + X \end{array}\right) \\
\left(\begin{array}{cc} -4a^2 & 0 \\ (aI+X)\ro & -(aI+X)^2 \end{array}\right) &=&
\left(\begin{array}{cc} -4a^2 & 0 \\ 2a\ro & -(aI+X)^2 \end{array}\right) \\
\phantom{\begin{array}{cc} 1&0\\ 0&1\end{array}}(aI+X)\ro &=& 2a\ro \\
X\ro &=& a\ro
\end{eqnarray*}
Thus for each $\al$, $\ro^\al$ is an eigenvector of the associated matrix $\Xa$ with eigenvalue $a_\al$.

By taking a linear combination of the $\Sa$'s we can perform a change of basis so that in the new basis $a_2 = a_3 = \cdots = a_f = 0$ and $a_1$ is either 1 or 0.  

As stated in \cite{rw}, there are at most $n$ linearly nilindependent matrices in $sp(2n,F)$.  This assertion follows from work in \cite{pwz}, and imposes a constraint on the number of  $\Sa$'s, namely $f \leq n+1$.   The dimension of $L$ is maximized when $f=n+1$.  In this case, $a_1=1$ and we can express the algebra with $X_1=0$ and $\{X_2, \ldots, X_f\}$ linearly nilindependent.  

We will now investigate in depth all cases where $a_1=1$.  We proceed by looking at the product of certain elements of the left-annihilator of $L$, $\Ann_\ell(L)$, with $S_1$ to place constraints on $\ro_1^\al$, $\ro_2^\al$, $r_{\al\be}$.
\begin{align*}
\begin{tabular}{ccc}
Identity	&	&	Constraints	\\
\hline
$[[\Sa,P_i]+[P_i,\Sa]\, ,S_1] = 0$	&	$\Rightarrow$	&	$\ro_1^\al = 0$\\
$[[\Sa,B_i]+[B_i,\Sa]\, ,S_1] = 0$	&	$\Rightarrow$	&	$\ro_2^\al = 0$\\
$[[\Sa,\Sb]+[\Sb,\Sa]\, ,S_1] = 0$	&	$\Rightarrow$	&	$r_{\al\be} = - r_{\be\al}$
\end{tabular}
\end{align*}
Thus $L_{\Sa} = -R_{\Sa}$ and $[\Sa,\Sb] = -[\Sb,\Sa]$, and so the algebra $L$ is a Lie algebra.  It should be noted that all Lie algebras with Heisenberg nilradical have been classified by Rubin and Winternitz in \cite{rw}.

Considering the Leibniz identity on the triple $\{S_1, \Sa, \Sb\}$ we get the following identity:
\begin{equation}\label{arar}
a_1r_{\al\be} = a_{\al}r_{1\be} - a_{\be}r_{1\al}
\end{equation}
By setting $\al, \be \neq 1$ in \eqref{arar} we arrive at $r_{\al\be}=0$. Setting $\al=\be=1$ gives $r_{11}=0$ and setting $\al\neq 1$ and $\be=1$ gives $r_{1\al}=0$. Hence, in the $a_1=1$ case we have that $r_{\al\be} = 0$ for all $\al$ and $\be$.  

We summarize in the following theorem.

\begin{theorem}\label{bigthm}
Every indecomposable solvable Leibniz algebra $L$ over $F=\C$ or $\R$ with nilradical $H(n)$ can be written in the basis $\{S_1, \ldots, S_f, P_1, \ldots, P_n,B_1,\ldots,B_n,H\}$ with multiplication relations \eqref{Heis}, $[\Sa,\Sb]=r_{\al\be} H$, and 
\begin{eqnarray*}\label{Ls2}
L_{\Sa} \left(
\begin{array}{c}
H\\
P\\
B
\end{array}
\right) &=& \left(
\begin{array}{ccc}
\phantom{-}2a_{\al}	&	0	&	0	\\
0	&	\phantom{-}a_\al I_n+\Aa	&	\Ca	\\
0	&	\Da	&	\phantom{-}a_\al I_n-A_\al^T
\end{array}
\right) \left(
\begin{array}{c}
H\\
P\\
B
\end{array}
\right)\\
\label{Rs2}
R_{\Sa} \left(
\begin{array}{c}
H\\
P\\
B
\end{array}
\right) &=& \left(
\begin{array}{ccc}
-2a_\al	&	0	&	0	\\
\ro_1^\al	&	-a_\al I_n-\Aa	&	-\Ca	\\
\ro_2^\al	&	-\Da	&	-a_\al I_n+\Aa^T 
\end{array}
\right) \left(
\begin{array}{c}
H\\
P\\
B
\end{array}
\right)
\end{eqnarray*}
\noindent for all $\al,\be =1,\ldots,f$.
The matrices $\Xa =  \left(\begin{array}{cc} \Aa&\Ca\\ \Da&-\Aa^T\end{array}\right) \in sp(2n,F)$ commute.  The constants $a_\al$ satisfy $a_1=1$ or $0$ and $a_2 = \cdots = a_f = 0$. If $a_1=1$, $\{X_2, \ldots, X_f\}$ is linearly nilindependent, $r_{\al\be}=0$, and $\ro_1^\al=\ro_2^\al=0$.  Hence, $L$ is a Lie algebra.  If $a_1=0$, $\{X_1, \ldots, X_f\}$ is linearly nilindependent, $r_{\al\be} \in F$, and $\ro^\al = (\ro_1^\al \,,  \ro_2^\al)^T$ is in the nullspace of $\Xa$.  The dimension of $L$ is $2n+1+f$ where $f \leq n+1$.
\end{theorem}

The isomorphism classes of the identification in the theorem are given as follows.  Two Leibniz algebras $L$ and $L^{\prime}$ both in the form of the theorem are isomorphic if and only if $L^{\prime}$ can be obtained from $L$ through one or more transformation that respects the form of the algebra given in the theorem.  Allowable transformations are: scaling the Heisenberg algebra
$$P_i^{\prime} = c P_i, \quad B_i^{\prime} = c B_i, \quad H_i^{\prime} = c^2 H_i,$$ where $c \neq 0 \in F$; symplectic transformations of $\{P_i,B_j\}$
$$\xi = \left(
\begin{array}{c}
P \\
B
\end{array}
\right), \quad \xi^{\prime}= G \xi,$$
with $GKG^{T}=K$ where $K = 
\left(
\begin{array}{cc}
0 & I_n\\
-I_n & 0
\end{array}
\right)$ as before, or equivalently $G \in Sp(2n,F)$; and linear combinations of $\{S_\alpha\}$ which respect the form of $L_{S_\alpha}$, $R_{S_\alpha}$, and $[S_\alpha, S_\beta]$.

In the following section, we will use the theorem to classify all solvable Leibniz algebras whose nilradical is the 3-dimensional Heisenberg algebra, $H(1)$.  The theorem can also be used to classify any solvable Leibniz algebra whose nilradical is Heisenberg.

\section{Extensions of $H(1)$}


By the restriction Theorem \ref{bigthm} imposes on the dimension of $L$, we need only consider 1- and 2-dimensional extensions of $H(1)$ (over the field $F=\C$ or $\R$).  First consider a 1-dimensional extension of the Heisenberg algebra $H(1)$ with basis $\{H, P, B\}$ by a new element $S$. 
Using the notation of the theorem, let $X=\left(\begin{array}{cc}A&C\\D&-A\end{array}\right)$, with $A,C,D\in F$. Then $L_S = \left(
\begin{array}{cc}
\phantom{-}2a	&	0	\\
0	&	\phantom{-}aI_2+X	
\end{array}
\right)$, $R_{S} = \left(
\begin{array}{cc}
-2a	&	0	\\
\ro	&	-aI_2-X	
\end{array}
\right)$,
and $[S,S]=rH$ where $\ro = (\ro_1\,,\ro_2)^T$ and $r,a,\ro_1,\ro_2\in F$.

If $a=1$, then $r = \ro_1 = \ro_2 = 0$, so $L$ is a Lie algebra. Rubin and Winternitz classified Lie algebras with Heisenberg nilradicals in \cite{rw}. These algebras are defined by the left-multiplication of $S$ on $H(1)$ via one of the following matrices,
\begin{equation}\label{H1a1C}
L_S = \left(\begin{array}{ccc}
	2	&		&		\\
		&	1+A	&	0	\\
		&	0	&	1-A	
\end{array}\right) \qquad,
\qquad
L_S = \left(\begin{array}{ccc}
	\ 2	&	\phantom{1+A}	&		\\
		&	1	&	1 \	\\
		&	0	&	1 \ 	
\end{array}\right)
\end{equation}
\begin{equation}\label{H1a1R}
L_S = \left(\begin{array}{ccc}
	\ 2	&	\phantom{1+A}	&		\\
		&	1	&	C \ 	\\
		&	-C	&	1 \ 	
\end{array}\right)
\end{equation}
where $A\ge 0$ and $C>0$.  The matrices in \eqref{H1a1C} represent the two distinct classes of algebras over $\C$.  The matrices in \eqref{H1a1C} and \eqref{H1a1R} represent the three distinct classes of algebras over $\R$.  Over $\C$, the algebras represented by the matrix on the left side of \eqref{H1a1C} and the matrix in \eqref{H1a1R} condense to just one case.

If $a=0$, then $S$ is linearly nilindependent.  Thus $S$ is not nilpotent, so $X \in sp(2,F)$ is not nilpotent.  With this, a simple calculation guarantees that the determinant of $X$ is nonzero.  Thus $X$ has a trivial nullspace.  Since $\ro$ is an eigenvector of $X$ with eigenvalue $a=0$, this gives that $\ro=0$.

If $r=0$, then $L$ is again a Lie algebra, so by \cite{rw} left-multiplication of $S$ on $H(1)$ is given by one of the following matrices.
\begin{equation}\label{H1a0C}
L_S = \left(\begin{array}{ccc}
	\ 0	&	\phantom{1+A}	&		\\
		&	1	&	0	\\
		&	0	&	-1 	
\end{array}\right)
\end{equation}
\begin{equation}\label{H1a0R}
L_S = \left(\begin{array}{ccc}
	\ 0	&	\phantom{1+A}	&	\phantom{-1}	\\
		&	0	&	1	\\
		&	-1	&	0 	
\end{array}\right)
\end{equation}
The matrix in \eqref{H1a0C} represents the only class of algebras over $\C$.  The matrices in \eqref{H1a0C} and \eqref{H1a0R} represent the two distinct classes of algebras over $\R$.  Over $\C$, these two classes condense to just one case.

If $r\neq 0$, then the change of basis $\tilde{S} = \dfrac{1}{\la}S$ where $\la^2 = \det(X)$ over $\C$ {\Large [}respectively $\la^2 = |\det(X)|$ over $\R${\Large ]} will make $\det(\tilde{X}) = -1$ {\Large [}resp. $\pm 1${\Large ]}.
Then $X$ is conjugate to $\left(\begin{array}{cc}1&0\\ 0&-1\end{array}\right)$ over $\C$ {\Large [}resp. $\left(\begin{array}{cc}1&0\\ 0&-1\end{array}\right)$ or $\left(\begin{array}{cc}0&1\\ -1&0\end{array}\right)$ over $\R${\Large ]}
so a suitable change of basis will put $X$ equal to this matrix {\Large [}resp. matrices{\Large ]}.
After these changes of basis we still have $r \neq 0$,
so we can scale the basis for $H(1)$ via $\tilde{P}=\mu P$, $\tilde{B}=\mu B$, and $\tilde{H}=\mu^2 H$, where $\mu^2=r$ {\Large [}resp. $\mu^2=|r|${\Large ]}.  This leaves the matrix $X$ unchanged but makes $[S,S]=rH=\dfrac{r}{\mu^2}\tilde{H}=\tilde{H}$ {\Large [}resp. $\pm \tilde{H}${\Large ]}.

Therefore the algebra will be defined by the value of $r$ and the left-multiplication of $S$ on $H(1)$.  Over $\C$, $r=1$ and left-multiplication is given by \eqref{H1a0C}.  Over $\R$, we have four cases where $r=\pm 1$ and left-multiplication can be given by either \eqref{H1a0C} or \eqref{H1a0R}.

It should be noted that in \cite{ck} Ca\~nete and Khudoyberdiyev have classified all 4-dimensional Leibniz algebras over $\C$, including the special case where the nilradical is Heisenberg.  
We have recovered their results with the algebras defined by \eqref{H1a1C}, and by \eqref{H1a0C} with $r=0$ or 1. We have further  classified all such algebras over $\R$, namely in \eqref{H1a1C} and \eqref{H1a1R}, and in \eqref{H1a0C} and \eqref{H1a0R} with $r = 0, 1,$ or $-1$.

Next, consider a 2-dimensional extension of $H(1)$ over $F$, with new basis elements $S_1$ and $S_2$, and corresponding matrices $X_1$ and $X_2$, using the notation in Theorem \ref{bigthm}. Specifically, recall that $a_2 = 0$.

If $a_1 = 1$, then $\rho_1^\al = \rho_2^\al = r_{\al\be} = 0$, for $\al,\be = 1$ or 2. This algebra is a Lie algebra. By \cite{rw}, the left-multiplications of these elements on $H(1)$ are given by the following matrices.
\begin{equation}\label{H2a1C}
L_{S_1} = 
\left(\begin{array}{ccc}
	\ 2	&	\phantom{1+A}	&		\\
		&	1	&	0\ 	\\
		&	0	&	1\ 	
\end{array}\right)\qquad ,\qquad
L_{S_2} = 
\left(\begin{array}{ccc}
	\ 0	&	\phantom{1+A}	&	\phantom{-1}	\\
		&	1	&	0	\\
		&	0	&	-1	
\end{array}\right)
\end{equation}
\begin{equation}\label{H2a1R}
L_{S_1} = 
\left(\begin{array}{ccc}
	\ 2	&	\phantom{1+A}	&		\\
		&	1	&	0\ 	\\
		&	0	&	1\ 	
\end{array}\right)\qquad ,\qquad
L_{S_2} = 
\left(\begin{array}{ccc}
	\ 0	&	\phantom{1+A}	&	\phantom{-1}	\\
		&	0	&	1	\\
		&	-1	&	0 	
\end{array}\right)
\end{equation}
The pair of matrices in \eqref{H2a1C} represents the only class of algebras over $\C$. The pairs of matrices in \eqref{H2a1C} and \eqref{H2a1R} represent the two distinct classes of algebras over $\R$. Over $\C$, these two classes condense to just one case.

%
If $a_1 = 0$, then by the same argument that we used in the 1-dimensional extension of $H(1)$ for $S_1$, we have that $\ro_1^1 = \rho_2^1 = 0$. Similarly, we argue that $\rho_1^2 = \rho_2^2 = 0$. Recalling that left-multiplication commutes we have that
\begin{eqnarray*}
	L_{S_1} L_{S_2} &=& L_{S_2} L_{S_1}	\\
	X_1 X_2 &=& X_2 X_1	\\
	\left(\begin{array}{cc} 
		A_1&C_1\\ D_1&-A_1
		\end{array}\right)
	\left(\begin{array}{cc} 
		A_2&C_2\\ D_2&-A_2
		\end{array}\right)
	&=&
	\left(\begin{array}{cc} 
		A_2&C_2\\ D_2&-A_2
		\end{array}\right)
	\left(\begin{array}{cc} 
		A_1&C_1\\ D_1&-A_1
		\end{array}\right).
\end{eqnarray*}
As a result, we observe that $\dfrac{A_1}{A_2} = \dfrac{C_1}{C_2} = \dfrac{D_1}{D_2}$. Hence, $X_1$ is a scalar multiple of $X_2$, and so they are not linearly nilindependent. 
Thus, there are no Leibniz algebras with $a_1=0$, so all 5-dimensional solvable Leibniz algebras with Heisenberg nilradical are Lie algebras and are given by \eqref{H2a1C} or \eqref{H2a1R}. This result is predicted by Theorem \ref{bigthm} since $f=n+1$.

We have now classified all indecomposable solvable Leibniz algebras over $\C$ and $\R$ whose nilradical is $H(1)$.

\vspace{12pt}

\noindent{\bf Acknowledgements.}  

The authors gratefully acknowledge the support of the faculty summer research grant from Spring Hill College, as well as the support of the Departments of Mathematics at Spring Hill College, the University of Texas at Tyler, and West Virginia University: Institute of Technology.


\begin{thebibliography}{ABCDE}


\bibitem{ao} Albeverio, Sh., A. Ayupov, B. Omirov.  On nilpotent and simple Leibniz algebras, \emph{Comm. Algebra}.  \textbf{33}  (2005), no.1,  159-172.

\bibitem{aor} Albeverio, S., B. A. Omirov, I. S. Rakhimov.  Classification of 4 dimensional nilpotent complex Leibniz algebras, \emph{Extracta Math.} \textbf{21} (2006), no. 3, 197-210.

\bibitem{ayupov} Ayupov, Sh. A., B. A. Omirov.  On Leibniz algebras, \emph{Algebra and Operator Theory}.  Proceedings of the Colloquium in Tashkent, Kluwer, (1998). 





\bibitem{barneslevi} Barnes, D.  On Levi's theorem for Leibniz algebras, \emph{Bull. Austral. Math. Soc.}  \textbf{86} (2012), no. 2, 184 - 185.

\bibitem{barnesleib} Barnes, D.  Some theorems on Leibniz algebras, \emph{Comm. Algebra}. \textbf{ 39} (2011), no. 7,  2463-2472.   

\bibitem{barnesengel} Barnes, D.  On Engel's Theorem for Leibniz Algebras, \emph{Comm. Algebra}. \textbf{40} (2012), no. 4, 1388-1389.


\bibitem{chelsie-allison} Batten-Ray, C., A. Hedges, E. Stitzinger.  Classifying several classes of Leibniz algebras, \emph{Algebr. Represent. Theory}.  To appear, arXiv:1301.6123.

\bibitem{jacobsonleib} Bosko, L., A. Hedges, J.T. Hird, N. Schwartz, K. Stagg.  Jacobson's refinement of Engel's theorem for Leibniz algebras, \emph{Involve}. \textbf{4}  (2011), no. 3, 293-296.

\bibitem{ck}Ca\~{n}ete, E. M., A. K. Khudoyberdiyev. The Classification of 4-Dimensional Leibniz Algebras. \emph{Linear Algebra Appl.} \textbf{439} 1, (2013) 273-288.

\bibitem{cs}Campoamor-Stursberg R. Solvable Lie algebras with an $\mathbb{N}$-graded nilradical of maximum nilpotency degree and their invariants, \emph{J. Phys A.} \textbf{43}  (2010), no. 14, 145202, 18 pp.

\bibitem{clok} Casas J. M., M. Ladra, B. A. Omirov, I. A. Karimjanov. Classification of solvable Leibniz algebras with null-filiform nilradical, \emph{Linear Multilinear Algebra}. \textbf{61}  (2013), no. 6, 758-774.

\bibitem{clok2} Casas J. M., M. Ladra, B. A. Omirov, I. A. Karimjanov. Classification of solvable Leibniz algebras with naturally graded filiform nil-radical, \emph{Linear Algebra Appl.} \textbf{438}  (2013), no. 7, 2973-3000.



\bibitem{loday} Loday, J.  Une version non commutative des algebres de Lie: les algebres de Leibniz, \emph{Enseign. Math. (2)}   \textbf{39}  (1993), no. 3-4, 269-293.

\bibitem{loday2} Loday, J.,  T. Pirashvili.  Universal enveloping algebras of Leibniz algebras and (co)homology, \emph{Math. Ann.} \textbf{296}  (1993), no. 1, 139-158.

\bibitem{mubar} Mubarakzianov, G. M. \emph{Some problems about solvable Lie algebras}. Izv. Vusov: Ser. Math. (1966), no. 1, 95-98 (Russian).

\bibitem{nw} Ndogmo, J. C., P. Winternitz.  Solvable Lie algebras with abelian nilradicals, \emph{J. Phys. A.} \textbf{27}  (1994), no. 8, 405-423.

\bibitem{omirov} Omirov, B.  Conjugacy of Cartan subalgebras of complex finite dimensional Leibniz algebras, \emph{J. Algebra}. \textbf{302}  (2006), no. 2, 887-896.

\bibitem{pwz} Patera, J., P. Winternitz, H. Zassenhaus. Maximal abelian subalgebra of real and complex symplectic Lie algebras, \emph{J. Math. Phys.}. \textbf{24}  (1983), 8, 1973-1985.

\bibitem{rw} Rubin, J. L., P. Winternitz.  Solvable Lie algebras with Heisenberg ideals, \emph{J. Phys. A: Math. Gen.} \textbf{26} (1993), no.  5, 1123-1138.

\bibitem{tw} Tremblay P., P. Winternitz.  Solvable Lie algebras with triangular nilradicals, \emph{J. Phys. A.} \textbf{31}  (1998), no. 2, 789-806.

\bibitem{wld} Wang, Y., J. Lin, S. Deng. Solvable Lie algebras with quasifiliform nilradicals, \emph{Comm. Algebra}. \textbf{36}  (2008), no. 11, 4052-4067.


\end{thebibliography}
\end{document}